\documentclass{amsart}\usepackage{amssymb}
\newtheorem{Lemma}{Lemma}[section]\newcommand{\bel}{\begin{Lemma}}\newcommand{\eel}{\end{Lemma}}
\newtheorem{Proposition}[Lemma]{Proposition}\newcommand{\bprop}{\begin{Proposition}}\newcommand{\eprop}{\end{Proposition}}
\newtheorem{Theorem}[Lemma]{Theorem}\newcommand{\bthe}{\begin{Theorem}}\newcommand{\ethe}{\end{Theorem}}
\newcommand{\bpr}{{\em Proof.\ }}\def\epr{$\bull$\\}
\newtheorem{Remark}[Lemma]{Remark}\newcommand{\beR}{\begin{Remark}\rm}\newcommand{\eeR}{\end{Remark}}
\newtheorem{Definition}[Lemma]{Definition}\newcommand{\bed}{\begin{Definition}}\newcommand{\eed}{\end{Definition}}
\newtheorem{Example}[Lemma]{Example}\newcommand{\bex}{\begin{Example}\rm}\newcommand{\eex}{\end{Example}}
\newtheorem{Corollary}[Lemma]{Corollary}\newcommand{\bcor}{\begin{Corollary}}\newcommand{\ecor}{\end{Corollary}}

\newcommand{\beq}{\begin{equation}}\newcommand{\eeq}{\end{equation}}
\newcommand{\bem}{\begin{displaymath}}\newcommand{\eem}{\end{displaymath}}
\newcommand{\beqa}{\begin{eqnarray}}\newcommand{\eeqa}{\end{eqnarray}}
\newcommand{\bee}{\begin{enumerate}}\newcommand{\eee}{\end{enumerate}}
\newcommand{\bei}{\begin{itemize}}\newcommand{\eei}{\end{itemize}}
\newcommand{\bet}{\begin{tabular}{cccccccc}}\newcommand{\eet}{\end{tabular}}
\newcommand{\bpm}{\begin{pmatrix}}\newcommand{\epm}{\end{pmatrix}}
\newcommand{\bM}{\begin{matrix}}\newcommand{\eM}{\end{matrix}}
\newcommand{\ber}{\begin{array}{l}}\newcommand{\eer}{\end{array}}

\def\bull{\vrule height .9ex width .9ex depth -.1ex }
\newcommand{\quotient}[2]{{\left.\raisebox{0.4ex}{$#1$}\!\!\middle/\!\!\raisebox{-0.4ex}{$#2$}\right.}}

\def\ra{\rightarrow}

\def\proj{\stackrel{\pi}{\ra}}

\def\li{~\\ $\bullet$ }

 \def\cK{{\mathcal{K}}}
\def\cO{\mathcal{O}}

\def\C{\mathbb{C}}

\def\k{k}

\def\Z{\mathbb{Z}}

  \def\tX{{\tilde{X}}}
\def\tY{\tilde{Y}}

\def\empty{\varnothing}
\def\suml{\sum\limits}

\def\prodl{\prod\limits}

\def\l{\ell}

\def\sset{\subset}

\def\?{{\bf ???}}\def\cf{C_{n,r}}
\newcommand{\bin}[2]{\binom{#1}{#2}}
\newcommand{\stir}[2]{\genfrac{\{}{\}}{0pt}{0}{#1}{#2}}

\def\uk{{\underline{k}}}\def\ux{{\underline{x}}}\def\uy{{\underline{y}}}

\title{T\MakeLowercase{he} `\MakeLowercase{corrected}
 D\MakeLowercase{urfee's inequality' for homogeneous complete intersections}}
\author{D\MakeLowercase{mitry} K\MakeLowercase{erner} \MakeLowercase{and} A\MakeLowercase{ndr\'as} N\MakeLowercase{\'emethi}}
\date{\today}

\address{Department of Mathematics, University of Toronto, 40 St. George Street, Toronto, Canada}
\email{dmitry.kerner@gmail.com}

\address{R\'enyi Institute of Mathematics\\Budapest\\Re\'altanoda u. 13\textendash15\\1053\\Hungary}
\email{nemethi@renyi.hu}

\thanks{A.N. is partially supported by OTKA grant 100796 of the
Hungarian Academy of Sciences.}
\thanks{Part of the work was done in Mathematische Forschungsinstitute
Oberwolfach, during D.K.'s stay as an OWL-fellow.}

\begin{document}\setcounter{secnumdepth}{6} \setcounter{tocdepth}{1}
\begin{abstract}
We address the conjecture of [Durfee1978], bounding the singularity genus $p_g$ by a multiple of the Milnor
number $\mu$ for an $n$--dimensional isolated complete intersection singularity. We show that the
original conjecture of Durfee, namely $(n +1)!\cdot  p_g \leq \mu$,
fails whenever the codimension $r$ is greater than one. Moreover,
we propose a new inequality $C_{n,r}\cdot  p_g \leq \mu$,
and  we verify it for homogeneous complete intersections.
In the homogeneous case
the inequality is guided by a `combinatorial inequality', that might have an independent interest.
\end{abstract}
\maketitle
\section{Introduction}

\subsection{} Let $(X,0)\sset(\C^N,0)$ be the analytic  germ of an $n$--dimensional  complex
 isolated complete intersection singularity (ICIS).
One of the most important goals  of the local singularity theory is
the clarification of the subtle connections between the two basic  numerical invariants,
the {\it Milnor number} $\mu$ and {\it singularity
(geometric) genus} $p_g$.

The surface case, $n=2$,  is already rather exotic and hard. In this case,
if $X_F$ is the Milnor fiber of $(X,0)$, and  $(\mu_+,\mu_0,\mu_-)$ are the
{\it Sylvester invariants} of the symmetric intersection form in the middle integral
homology  $H_2(X_F,\Z)$, then $2p_g=\mu_0+\mu_+$ \cite{Durfee1978},
while, obviously, $\mu=\mu_++\mu_0+\mu_-$. Hence, numerical relations between $\mu$ and $p_g$
can be rewritten in terms of the Sylvester invariants too. In topology one  also uses the
{\it signature} $\sigma:=\mu_+-\mu_-$ as well. In fact, for compact complex surfaces,
the Euler number, Todd genus and the signature are the most important index--theoretical
numerical invariants; their local analogs are the above integers $\mu$, $p_g$ and $\sigma$.
For more about these invariants see the monographs \cite{Milnor-book,AGLV,Looijenga-book}  
or the articles \cite{Laufer1977,Looijenga1986,Saito1981}.
For various formulae regarding the Milnor number of weighted homogeneous complete intersections see
\cite{Greuel1975,Greuel-Hamm1978,Hamm1986,Hamm2011} and for the geometric genus see
\cite{Khovanskii1978,Morales1985}.

Examples show that for  a local surface singularity $\mu_-$ should be `large' compared with the
other Sylvester invariants, or equivalently, $p_g$ small with respect to $\mu$, or, $\sigma$ rather negative.

This was formulated more precisely in {\bf Durfee's Conjectures}  \cite{Durfee1978} as follows:

\vspace{2mm}

\noindent {\bf (A)} \
 {\em Strong inequality}: if $(X,0)$ is an isolated complete intersection surface singularity, then $6p_g\leq \mu$.\\
 {\bf (B)} \ {\em Weak inequality}: if $(X,0)$ is a normal surface singularity (not necessarily ICIS)
 which admits a smoothing with Milnor number (second Betti number of the fiber) $\mu$,
then $4p_g\leq \mu+\mu_0 $, or equivalently, $\sigma\leq 0$.\\
{\bf (C)} \ {\em Semicontinuity of $\sigma$}: if $\{(X_t,0)\}_{t\in (\C,0)}$ is a family of
isolated surface singularities then $\sigma(X_{t=0})\leq \sigma(X_{t\not=0})$.

\vspace{2mm}

Almost immediately a counterexample to the {\it weak inequality }  was given in
\cite[page 240]{Wahl 1981} providing a
normal surface singularity (not ICIS) with $\mu=3$, $\mu_0=0$  and $p_g=1$.

A counterexample to the semicontinuity of the signature was found much later, in \cite{Kerner-Nemethi2009}.

On the other hand, the {\it strong inequality},  valid  for an ICIS,
 was believed to be true  and was verified for many
particular {\it hypersurface singularities $(X,0)$ in $(\C^3,0)$}:

\vspace{2mm}

\cite{Tomari93} proved $8p_g<\mu$ for $(X,0)$ of multiplicity 2,

\cite{Ashikaga92} proved $6p_g\le\mu-2$ for $(X,0)$ of multiplicity 3,

\cite{Xu-Yau1993} proved $6p_g\le\mu-mult(X,0)+1$ for quasi-homogeneous singularities,

\cite{Nemethi98,Nemethi99} proved $6p_g\le\mu$ for suspension type singularities $\{g(x,y)+z^k=0\}\sset(\C^3,0)$,

\cite{Melle00} proved $6p_g\le\mu$ for absolutely isolated singularities.

\vspace{2mm}

Moreover, for arbitrary $n\geq 2$, \cite{Yau-Zhang2006} proved the inequality $(n+1)!p_g\leq \mu$
for {\it isolated weighted-homogeneous hypersurface singularities} in $(\C^{n+1},0)$. The natural expectation
 was that the same inequality  holds for any ICIS of any dimension $n$ and any
 codimension $r:=N-n$.

\subsection{} This paper is the continuation of \cite{Kerner-Nemethi.a}. 
 The main results of the present article are the  following:

\vspace{2mm}

\noindent {\bf I.} \  For homogeneous ICIS of multidegree $(p_1,\ldots, p_r)$
we provide new formulae for $\mu$ and $p_g$; their special form allows us to compare them.\\
{\bf II.} \  Using these formulae one sees rather easily that $(n+1)!\cdot p_g\leq \mu$ is not true
whenever $r\geq 2$ and $n\geq 2$, already for $p_1=\cdots=p_r$ sufficiently large (that is, the strong   inequality
fails even asymptotically). \\
{\bf III.} \ We propose a new set of conjectured inequalities with new bounds for $\mu/p_g$.
For any $n\geq 1$ and $r\geq 1$
consider the Stirling number of the second kind,   cf. \cite[\S24.1.4]{Abramowitz-Stegun}, and the
coefficient $\cf$ defined by
$$\stir{n+r}{r}:=\frac{1}{r!}\suml^r_{j=0}(-1)^j\bin{r}{j}(r-j)^{n+r},$$
\beq\label{eq:1} \cf:=\frac{\bin{n+r-1}{n}(n+r)!}{\stir{n+r}{r}r!}=
\frac{|\cK_{n,r}|}{\suml_{\uk\in\cK_{n,r}}\prod^r_{i=1}\frac{1}{(k_i+1)!}},\eeq
where $|\cK_{n,r}|=\binom{n+r-1}{n}$   is the cardinality of the set
\beq
\cK_{n,r}:=\{\uk=(k_1,\ldots, k_r)\, :\,  k_i\geq 0 \ \mbox{for all $i$, and } \ \sum _ik_i=n\}.
\eeq
The second equality of (\ref{eq:1}) follows from \cite[pages 176-178]{Jordan1965}.
One also shows, cf. Corollary \ref{ss:ineq},  that
\beq\label{eq:C<}
C_{n,1}> C_{n,2}> \cdots > C_{n,r} >\cdots > \lim_{r\to\infty}\cf.
\eeq
E.g., for small $r$ and for $r\to\infty$ one gets
\beq
C_{n,1}=(n+1)! \ \ \ \ \ \ C_{n,2}=\frac{(n+2)!(n+1)}{2^{n+2}-2} \ \ \ \ \ \ \lim_{r\to\infty}\cf=2^n.
\eeq
The limit can be  computed using the asymptotical growth  of Stirling numbers of the second kind,
\cite[\S24.1.4]{Abramowitz-Stegun}:
$\stir{n+r}{r}\sim\frac{r^{2n}}{2^nn!}$. This gives: $\cf\sim 2^n\frac{(n+r-1)!(n+r)!}{(r-1)!r!r^{2n}}$
with limit $2^n$ as $r\to\infty$.

\vspace{2mm}

\noindent
{\bf Conjecture.} {\em Let $(X,0)\sset(\C^N,0)$ be an ICIS of dimension $n$ and of codimension $r=N-n$. Then
\li for $n=2$ and $r=1$ one has  \ $6p_g\leq\mu$,
\li for $n=2$ and arbitrary $r$ one has  \ $4p_g<\mu $ (see last section too),
\li for $n\ge3$ and arbitrary  $r$ one has  \ $\cf\cdot  p_g\leq \mu $.}

\vspace{2mm}

\noindent {\bf IV.} \  We show that the third proposed inequality is asymptotically sharp,
 i.e. for any fixed $n$ and $r$ there exists a sequence of isolated
complete intersections for which the ratio $\frac{\mu}{p_g}$ tends to $\cf$. \\
{\bf V.} \ We support the above conjecture by its proof for any homogeneous ICIS with any
multidegree $(p_1,\ldots, p_r)$.

\vspace{2mm}

Note that $\cf\cdot  p_g\leq \mu $ automatically implies
$$\mu\geq \min_{\uk\in\cK_{n,r}}\{\prod_{i=1}^r (k_i+1)!\}\cdot p_g\geq 2^n\cdot p_g.$$

Some more comments are in order.

\li The general definition of the singularity genus is the
following. Let $(X,0)$ be a reduced isolated singularity of
dimension $n$, and let $\tX\proj(X,0)$ be a resolution. Let
$\cO_{(X,0)}$ be the local ring of the singulaity $(X,0)$, and let
$\cO_\tX$ be the structure sheaf on $\tX$.  Then  $\pi_*\cO_\tX$
is the normalization of $\cO_{(X,0)}$ and (cf.
\cite[\S4.1]{Looijenga1986}):
\[(-1)^np_g(X,0):=\sum_{i\ge1}(-1)^{i-1}h^i(\cO_\tX)-dim_\C\quotient{\pi_*\cO_\tX}{\cO_{(X,0)}}.\]
$\bullet$ For $n=1$ the preimage of the singular point is of
dimension zero, hence the higher cohomologies vanish:
$h^i(\cO_\tX)=0$ for all $i\geq 1$. Thus the singularity genus
coincides with the classical {\em delta invariant} $\delta$.  It
satisfies $2\delta=\mu+b(X,0)-1$, where $b(X,0)$ is the number of
locally irreducible components,
 \cite[Theorem 10.5]{Milnor-book}.
Hence $\mu\le 2\delta$, and the analogue of the  strong inequality fails (note that $C_{1,r}=\frac{1}{2}$). Moreover,
for $n=2$ too, the inequality   $\mu(X,0)\ge\cf\cdot  p_g(X,0)$, in general,  is not satisfied,
and the asymptotically sharp inequality  of this form is impossible.
For more comments regarding $n=2$ see section \ref{s:4}.

These facts also show that, in order to prove the conjectured inequalities,
a `naive induction' over $n$  is impossible.
\li  Homogeneous singularities (considered in {\bf V}) are
rather  particular ones.   Nevertheless, it turns out
that they are the building blocks for many other singularity types.
In the forthcoming paper we use the statement of {\bf V}
to prove the above conjectured   inequalities from {\bf III}  for
more general families of complete intersections (e.g., for absolutely isolated singularities),
cf. \cite{Kerner-Nemethi.c}.
\li Even in this particular case of homogeneous germs
the proof involves a non-trivial combinatorial inequality, which `guides'
the inequality $\cf \cdot p_g\leq \mu$, see \S1.3 and \S4.
Although its proof in its current form is relatively short, we believe it is far from being straightforward.

\subsection{}
Fix the integers $n\geq 1$, $r\geq 1$ and $\ell\geq 0$. We define $\uy_{\uk,\ell}:=\prod_i\frac{1}{(k_i+\ell)!}$
and for
free positive variables $x_1,\ldots, x_r$ we set  $\ux^{\uk}=x_1^{k_1}\cdots x_r^{k_r}$.
By convention $0!=1$.

\vspace{2mm}

\noindent {\bf Combinatorial Inequality.} \ {\it
For any $n,r$ and $\ell$ as above, and for any positive $x_1,\cdots, x_r$ one has:}
\beq (I_{n,r,\ell}):\hspace{1cm}
\frac{1}{|\cK_{n,r}|}\cdot \sum_{\uk\in\cK_{n,r}} \uy_{\uk,\ell}\cdot \sum_{\uk\in\cK_{n,r}} \ux^{\uk}\ge
\sum_{\uk\in\cK_{n,r}} \uy_{\uk,\ell}\cdot \ux^{\uk}.
\eeq

\vspace{2mm}

For our application  we only need the $\ell=1$ case, nevertheless its proof uses all the $\ell$ values
(the induction over $n$  involves larger $\ell$'s for smaller $n$'s).

\subsection{}
We wish to thank G.-M. Greuel, H. Hamm,  A. Khovanskii, M. Leyenson, L. Lov\'asz, P. Milman, E. Shustin for advices and
 important discussions.

\section{The $\mu$ and $p_g$ formulae}

Let $(X,0)=\{f_1=\cdots=f_r=0\}$ be a homogeneous ICIS of multidegree $(p_1,\ldots, p_r)$,
that is $deg(f_i)=p_i$.
\bprop\label{Thm.Milnor.Genus.For.ICIS}
(1) \ \ \ $\mu=\Big(\prodl^r_{i=1}p_i\Big)\suml^n_{j=0}(-1)^j
\bigg(\suml_{\uk\in\cK_{n-j,r}}
\prod_{i=1}^r(p_i-1)^{k_i}\bigg)-(-1)^n$,
\\ \hspace*{3cm}(2) \ \ \
$p_g=\suml_{\uk\in\cK_{n,r}}\prodl^r_{i=1}\bin{p_i}{k_i+1}$.
\eprop
The reader is invited to consult
 \cite{Greuel1975,Greuel-Hamm1978,Hamm1986,Hamm2011} for $\mu$ and
\cite{Khovanskii1978,Morales1985,Hamm2011} for $p_g$ of a weighted homogeneous ICIS.
These  formulae usually are rather different
 than ours considered above; nevertheless, both formulae (1) and (2) can be derived from expressions
 already present in the literature (though we have found them in a different way).

\bpr (1) We determine the Euler characteristic $\chi=(-1)^n\mu +1$ of the Milnor fiber.
For a power series $Z:=\sum_{i\geq 0}a_ix^i$ we write $[Z]_n$ for the coefficient $a_n$ of $x^n$.
By 3.7(c) of  \cite{Greuel-Hamm1978}
$$\chi=\prod_{i=1}^rp_i\cdot \Big[\frac{(1+x)^N}{\prod_i(1+p_ix)}\Big]_n.$$
Rewrite $1+p_ix$ as $(1+x)(1-\frac{(1-p_i)x}{1+x})$, hence
$$\Big[\frac{(1+x)^N}{\prod_i(1+p_ix)}\Big]_n=
\Big[(1+x)^n\cdot \prod_i\sum_{k_i\geq 0}\Big(\frac{(1-p_i)x}{1+x}\Big)^{k_i}\Big]_n=$$
$$\Big[ \sum _{k_1\geq 0,..,k_r\geq 0}x^{\sum k_i}(1+x)^{n-\sum k_i}\prod_i(1-p_i)^{k_i}\Big]_n=
 \sum _{k_1\geq 0,..,k_r\geq 0\atop \sum k_i\leq n}\prod_i(1-p_i)^{k_i}.
$$
(2) By \cite[Theorem 2.4]{Morales1985} (and computation of the number of lattice points under the
`homogeneous Newton diagram')
\beq\label{pg}
p_g=\bin{\sum_kp_k}{N}-\sum_{1\leq i\leq r}\bin{(\sum_kp_k)-p_i}{N}+
\sum_{1\leq i<j\leq r}\bin{(\sum_kp_k)-p_{i}-p_{j}}{N}-\ldots
\eeq
Using the  Taylor expansion $\frac{1}{(1-z)^{N+1}}=\sum_{l\geq N}\bin{l}{N}z^{l-N}$,
the right hand side of (\ref{pg})  is $\Big[\frac{\prod_i (1-z^{p_i})}{(1-z)^{N+1}}\Big]_{\sum_kp_k-N}$.
Thus
$$
p_g=rez_{z=0}F(z), \ \ \mbox{where} \ \
F(z):=\frac{\prod_i (1-z^{p_i})}{z^{\sum p_i-N+1}(1-z)^{N+1}}.
$$
Note that $rez_{z=\infty}F(z)=0$ since  $F(1/z)/z^2$ is regular at zero.
Hence,  $p_g=-rez_{z=1}F(z)$ too.
By the change of variables $z\mapsto 1/z$, this last expression transforms into
$$
p_g=rez_{z=1}\frac{\prod_{i=1}^r (z^{p_i}-1)}{(z-1)^{N+1}}.
$$
Since $z^p-1=\sum_{k\geq 0}\binom{p}{k+1}(z-1)^{k+1}$, we obtain
$$
p_g
=rez_{z=1}\frac{1}{(z-1)^{n+1}}\cdot \prod_{i=1}^r \sum_{k_i\geq 0}\binom{p_i}{k_i+1}(z-1)^{k_i}=
\sum_{\uk\in\cK_{n,r}}\binom{p_1}{k_1+1}\cdots \binom {p_r}{k_r+1}.\ \
\ \ \ \bull
$$
Consider the particular case $p_1=\cdots=p_r=p$. Then the above formulae
read as
\beq\label{eq:pp}\ber
\mu=(-1)^n\Big(p^r\sum_{j=0}^n(1-p)^j\bin{j+r-1}{j}-1\Big), \ \mbox{(see also \cite[3.10(b)]{Greuel-Hamm1978})},
\\\\p_g=\sum_{\uk\in\cK_{n,r}}\prod^r_{i=1}\bin{p}{k_i+1}.
\eer\eeq
Therefore, for  $p$ large,  $\mu$ and $p_g$ asymptotically behave as follows:
 $$ \mu= p^N\bin{N-1}{n}+O(p^{N-1}) \ \ \ \mbox{and } \ \ \
p_g=p^N\cdot\suml_{\uk\in\cK_{n,r}}\prodl^r_{i=1}\frac{1}{(k_i+1)!}+O(p^{N-1}).$$
Thus, asymptotically, $\frac{\mu}{p_g}=C_{n,r}+O(\frac{1}{p})$.
Note that $C_{n,r}<C_{n,1}=(n+1)!$ for any $r\geq 2$, hence {\it
the strong Durfee's inequality  is violated  for any $p$
sufficiently large whenever $r\geq 2$}.

The inequality $C_{n,r}<(n+1)!$ is the consequence of  Corollary
\ref{ss:ineq}, but it follows from the next elementary observation
as well: $\frac{1}{(n+1)!}$ is the smallest of all elements of
type $\prod_i\frac{1}{(k_i+1)!}$, for $\uk\in\cK_{n,r}$. Indeed,
note that $(k_1+1)!(k_2+1)!\le(k_1+k_2+1)!$, because
 $\Big(\underbrace{2\times\cdots\times(k_1+1)}_{k_1}\Big)\times\Big(\underbrace{2\times\cdots\times(k_2+1)}_{k_2}\Big)
 \le\Big(\underbrace{2\times\cdots\times(k_1+k_2+1)}_{k_1+k_2}$\Big), thus $(n+1)!\ge\prod_i(k_i+1)!$.

\section{The inequality $\cf\cdot  p_g\leq \mu $ in the homogeneous case}\label{s:2}

Note that if $r=1$, then $ C_{n,1} \cdot p_g=(n+1)!\cdot p_g\leq
\mu$ for any isolated homogeneous germs and for any $n\geq 2$.
This follows from  $(n+1)!p_g=(n+1)!\binom{p}{n+1}\leq
(p-1)^{n+1}=\mu $, where $p$ is the degree of the germ. If $r>1$
then the $n=2$ case is rather special, and it will be discussed in
the last section. Hence, we start the case $n\geq 3$. We prove:
\bthe\label{the:In} Assume that  $(X,0)\sset(\C^{n+r},0)$ is a
homogeneous ICIS  of dimension $n>2$, codimension $r$ and
  multidegree $(p_1,\ldots, p_r)$.  Then $ \cf \cdot p_g\leq \mu$.
\ethe
In the next discussions  we assume $p_i\ge2$ for all $i$; if $p_i=1$ for some $i$ then one can reduce
the setup  to the  smaller $(r-1)$--codimensional case.
In addition, as all the formulas are symmetric in $\{p_i\}_i$, we will
sometimes assume that {\it $p_r$ is largest among all the $p_i$'s}.
Hence we wish to prove:
\beq\label{Eq.:3}
LHS:=(\prodl^r_{i=1} p_i) \suml^n_{j=0}(-1)^j\bigg(\suml_{\uk\in\cK_{n-j,r}}
\prod_{i=1}^r(p_i-1)^{k_i}\bigg)-(-1)^n\ge\cf\cdot
(\prodl^r_{i=1}p_i)\cdot \suml_{\uk\in\cK_{n,r}}\prod^r_{i=1}\bin{p_i-1}{k_i}\frac{1}{k_i+1}=:RHS.
\eeq
The proof consists of several steps.
\bthe
(1) The inequality (\ref{Eq.:3}) is the consequence of the next inequality:
$$\sum_{\uk\in\cK_{n,r}}\prod_{i=1}^r(p_i-1)^{k_i}\ge\cf\Bigg(
\sum_{\uk\in\cK_{n,r}}\prod^r_{i=1}\bin{p_i-1}{k_i}\frac{1}{k_i+1}+
\sum_{\substack{\uk\in\cK_{n,r}\\k_r>0}}\bin{p_r-2}{k_r-1}\frac{1}{k_r(k_r+1)}
\prod^{r-1}_{i=1}\bin{p_i-1}{k_i}\frac{1}{k_i+1}\Bigg).
$$
This inequality is the consequence of two further inequalities, listed in part (2) and (3):

\vspace{2mm}

\noindent (2) For $n\ge3$ the following inequality holds:
$$\sum_{\uk\in\cK_{n,r}}\prod_{i=1}^r\frac{(p_i-1)^{k_i}}{(k_i+1)!}\ge
\sum_{\uk\in\cK_{n,r}}\prod^r_{i=1}\bin{p_i-1}{k_i}\frac{1}{k_i+1}+
\sum_{\substack{\uk\in\cK_{n,r}\\k_r>0}}\bin{p_r-2}{k_r-1}\frac{1}{k_r(k_r+1)}
\prod^{r-1}_{i=1}\bin{p_i-1}{k_i}\frac{1}{k_i+1}.
$$
(3) For any $n\geq 2$ and $r\geq 1$  the following inequality holds:
$$\sum_{\uk\in\cK_{n,r}}\prod_{i=1}^r(p_i-1)^{k_i}\ge\cf\sum_{\uk\in\cK_{n,r}}\,
\prod_{i=1}^r\frac{(p_i-1)^{k_i}}{(k_i+1)!}.$$
\ethe
Here the third statement (3) is the most complicated one, it follows from the general
Combinatorial Inequality from the Introduction, and it is proved separately in the next section.

\vspace{2mm}

\noindent \bpr {\bf (1)} \  Expand the {\it LHS} of (\ref{Eq.:3})
in terms with decreasing $r$: \beq\label{Eq.:9}
LHS=(p_r-1)(\prodl^{r-1}_{i=1}p_i)C_n(p_1,\ldots,p_{r})+
(p_{r-1}-1)(\prodl^{r-2}_{i=1}p_i)C_n(p_1,\ldots,p_{r-1})+\cdots+(p_1-1)(p_1-1)^n,
\eeq where
$C_n(p_1,\ldots,p_{s}):=\suml_{\uk\in\cK_{n,s}}\,\prodl^s_{i=1}(p_i-1)^{k_i}$,
for $1\leq s\leq r$. To prove this formula we observe:
\[\ber
\suml^n_{j=0}(-1)^j\bigg(\suml_{\uk\in\cK_{n-j,r}}\prod_{i=1}^r(p_i-1)^{k_i}\bigg)=
(-1)^n\suml_{\substack{\suml^r_{j=1}k_j\le
n\\k_j\ge0}}\prodl^r_{i=1}(1-p_i)^{k_i}=\\=
(-1)^n\suml_{\substack{\suml^{r-1}_{j=1}k_j\le
n\\k_j\ge0}}\Big(\suml^{n-\suml^{r-1}_{i=1}k_i}_{k_r=0}(1-p_r)^{k_r}\Big)
\prodl^{r-1}_{i=1}(1-p_i)^{k_i} =
(-1)^n\suml_{\substack{\suml^{r-1}_{j=1}k_j\le
n\\k_j\ge0}}\Big(\frac{(1-p_r)^{n+1-\sum^{r-1}_{i=1}k_i}-1}{1-p_r-1}\Big)
\prodl^{r-1}_{i=1}(1-p_i)^{k_i}
\eer\]
Thus
\[
(\prodl^r_{i=1}p_i)(-1)^n\suml_{\substack{\suml^r_{j=1}k_j\le n\\k_j\ge0}}\prodl^r_{i=1}(1-p_i)^{k_i}=
(p_r-1)(\prodl^{r-1}_{i=1}p_i)C_n(p_1,\ldots,p_{r})+
(\prodl^{r-1}_{i=1}p_i)(-1)^n\suml_{\substack{\suml^{r-1}_{j=1}k_j\le n\\k_j\ge0}}\prodl^{r-1}_{i=1}(1-p_i)^{k_i}
\]
Iterating this gives equation (\ref{Eq.:9}).

It is natural to expand the right hand side of (\ref{Eq.:3}) similarly. For this, define
$D_n(p_1,\ldots,p_s):=
\suml_{\uk\in\cK_{n,s}}\prod^s_{i=1}\bin{p_i-1}{k_i}\frac{1}{k_i+1}$.
E.g.,  $D_1(p_1)=\bin{p_1-1}{n}\frac{1}{n+1}$.
Then, we write  $RHS/\cf$ as
$$
\Big((\prod^r_{i=1}p_i)D_n(p_1,\ldots,p_r)-(\prod^{r-1}_{i=1}p_i)D_n(p_1,\ldots,p_{r-1})\Big)+
\Big((\prod^{r-1}_{i=1}p_i)D_n(p_1,\ldots,p_{r-1})-(\prod^{r-2}_{i=1}p_i)D_n(p_1,\ldots,p_{r-2})\Big)
+\cdots+p_1 D_1(p_1).
$$
Thus, it is enough to prove the inequality for each pair of terms in these expansions, namely:
$$
 (p_s-1)(\prodl^{s-1}_{i=1}p_i)C_n(p_1,\ldots,p_s)\ge\cf
\Bigg((\prod^s_{i=1}p_i)D_n(p_1,\ldots,p_s)-(\prod^{s-1}_{i=1}p_i)D_n(p_1,\ldots,p_{s-1})\Bigg), \
\mbox{for all $1\leq s\leq r$}.
$$
Since  $C_{n,r}\leq C_{n,s}$, cf (\ref{eq:C<}), it is enough to prove the last inequality with coefficient
$C_{n,s}$ instead of $C_{n,r}$, or equivalently,
$$(p_r-1)(\prodl^{r-1}_{i=1}p_i)C_n(p_1,\ldots,p_r)\ge\cf
\Bigg((\prod^r_{i=1}p_i)D_n(p_1,\ldots,p_r)-(\prod^{r-1}_{i=1}p_i)D_n(p_1,\ldots,p_{r-1})\Bigg)
\ \mbox{for all $r\ge1,\ n\ge1$}.
$$
Further, split the right hand side of this last inequality into two parts:
$$
(p_r-1)(\prod^{r-1}_{i=1}p_i)D_n(p_1,\ldots,p_r)+
(\prod^{r-1}_{i=1}p_i)\Big(D_n(p_1,\ldots,p_r)-D_n(p_1,\ldots,p_{r-1})\Big),
$$
and rewrite $D_n(p_1,\ldots,p_r)-D_n(p_1,\ldots,p_{r-1})$ into
$$
\suml_{\substack{\uk\in\cK_{n,r}\\k_r>0}}\prod^r_{i=1}\bin{p_i-1}{k_i}\frac{1}{k_i+1}=
(p_r-1)\sum_{\substack{\uk\in\cK_{n,r}\\k_r>0}}\bin{p_r-2}{k_r-1}\frac{1}{k_r(k_r+1)}
\prod^{r-1}_{i=1}\bin{p_i-1}{k_i}\frac{1}{k_i+1}.
$$
This provides precisely the expression of part (1).

\vspace{2mm}

\noindent {\bf (2)} \ We start to compare individually the particular summands indexed by
$\uk\in\cK_{n,r}$ of both sides. First, we consider some $\uk\in\cK_{n,r}$ with  $k_r>1$.
Then the  corresponding individual inequality is true. Indeed,
$$\ber
\frac{\prodl^r_{i=1}\frac{(p_i-1)^{k_i}}{k_i!}-\prodl^r_{i=1}\bin{p_i-1}{k_i}}{\prodl^r_{i=1}(k_i+1)}-
\bin{p_r-2}{k_r-1}\frac{1}{k_r(k_r+1)}\prodl^{r-1}_{i=1}\bin{p_i-1}{k_i}\frac{1}{k_i+1}\ge
\frac{\frac{(p_r-1)^{k_r-1}}{(k_r-1)!}-\bin{p_r-2}{k_r-1}}{k_r(k_r+1)}\prodl^{r-1}_{i=1}\bin{p_i-1}{k_i}\frac{1}{k_i+1}\ge0.
\eer$$
Here and in the sequel we constantly use
$(p_i-1)^{k_i}-(p_i-1)\cdot(p_i-2)\cdots(p_i-k_i)\ge(p_i-1)^{k_i-1}$, valid for $k_i>1$.

Next, we assume that $\uk$ has the following properties: $k_r=1$, but $k_j>1$ for some $j$. Then again
$$\ber
\frac{p_r-1}{2}\cdot\frac{\prodl^{r-1}_{i=1}\frac{(p_i-1)^{k_i}}{k_i!}-\prodl^{r-1}_{i=1}\bin{p_i-1}{k_i}}
{\prod^r_{i=1}(k_i+1)}-\frac{1}{2}\prodl^{r-1}_{i=1}\bin{p_i-1}{k_i}\frac{1}{k_i+1}\ge
\frac{p_r-1}{2}\cdot\frac{(p_j-1)^{k_j-1}}{(k_j+1)!}\prodl_{\substack{1\leq i<r\\i\neq j}}
\frac{(p_i-1)^{k_i}}{(k_i+1)!}-\frac{1}{2}\prodl^{r-1}_{i=1}\bin{p_i-1}{k_i}\frac{1}{k_i+1}\ge0.
\eer$$
Here we used the initial assumption that $p_r\ge p_i$ for any $i$.

Now, we assume that $k_r=1$ and  $k_i\le1$ for all $i$.
This can happen only for $n\le r$.  Let  $i_1,\dots,i_{n-1}$ be those indices different than $r$
for which $k_i=1$, that is,  $k_r=k_{i_1}=\cdots=k_{i_{n-1}}=1$ and all the other $k_i$'s are zero.

In this case $\prod_{i=1}^r\frac{(p_i-1)^{k_i}}{(k_i+1)!}=\prod^r_{i=1}\bin{p_i-1}{k_i}\frac{1}{k_i+1}$,
hence the individual inequality corresponding to $\uk$ fails. Therefore, we will group this
term by some other terms with $k_r=0$.  More precisely, we will group  this $\uk$ together with
terms which correspond to those  $\uk$'s which satisfy
$k_r=0$, \ $k_{i_j}=2$ for exactly one $j\in \{1,\ldots, n-1\}$,
$k_{i_{l}}=1$ if $l\in \{1,\ldots, n-1\}\setminus \{j\}$,  and all the other $k_i$'s are zero.

Note that if $k_r=0$ then there  is no contribution from the sum
$\suml_{\substack{\uk\in\cK_{n,r}, \ k_r>0}}$.
Therefore, the $n$ individual inequalities corresponding to the above $\uk$'s altogether provide
$$
\suml^{n-1}_{l=1}\frac{(p_{i_l}-1)-(p_{i_l}-2)}{3}\cdot\prod^{n-1}_{l=1}\frac{p_{i_l}-1}{2}
-\frac{1}{2}\prod^{n-1}_{l=1}\frac{p_{i_l}-1}{2}=
\frac{2n-5}{6}\cdot\prod^{n-1}_{l=1}\frac{p_{i_l}-1}{2}.
$$
For $n\ge3$ this is positive, hence the statement.

Any other remaining $\uk\in\cK_{n,r}$ can again be treated individually: for all of them $k_r=0$ and
$\frac{(p-1)^k}{(k+1)!}\geq \binom{p-1}{k}\frac{1}{k+1}$.
\epr

\section{Proof of the Combinatorial Inequality}
We use the notations of \S1.3 and introduce more objects.
 We will consider the following partition of
$\cK_{n,r}$: for any $s\in \{0, \ldots, r-1\}$ we  define
\beq  \cK_{n,r}^s:=\{\uk\in \cK_{n,r}\, :\, |\{i:\, k_i=0\}|=s\}.
\eeq
Note that for $s<r-n$ one has
$\cK^{s}_{n,r}=\empty$.
Corresponding to these sets, we consider the {\it arithmetic mean}
$X_{n,r}$ and  $X^s_{n,r}$ of the elements $\ux^{\uk}$  indexed by the sets
  $\cK_{n,r}$ and $\cK_{n,r}^s$ respectively.
In parallel,  $Y_{n,r,\ell}$ and  $Y^s_{n,r,\ell}$ denote the
arithmetic mean of elements $\uy_{\uk,\ell}$ indexed by the same sets
$\cK_{n,r}$  and $\cK_{n,r}^s$ respectively.

\bthe\label{the:2} With the above notation one has

(a) \ $X^0_{n,r}\leq X^1_{n,r}\leq \cdots \leq X^{r-1}_{n,r}$

(b) \ $Y^0_{n,r,\ell} > Y^1_{n,r,\ell}> \cdots > Y^{r-1}_{n,r,\ell}$

(c) \ (a) and (b) imply the Combinatorial Inequality from the Introduction.
\ethe
\bpr {\bf (a)} \
By definition,  $X^{s}_{n,r}$ is the arithmetic mean of $\bin{r}{s}\bin{n-1}{r-s-1}$ monomials. Let
$\tX^{s}_{n,r}(x_1,\ldots,x_r)$ be the sum of these monomials, that is,
 $\tX^{s}_{n,r}(x_1,\ldots,x_r)=X^{s}_{n,r}\cdot\bin{r}{s}\bin{n-1}{r-s-1}$.

{\bf Step 1.}
We show that all the inequalities can be deduced  from  the first
one: $X^1_{n,r}\ge X^0_{n,r}$. Indeed,
by definition,  $\tX^{s}_{n,r}(x_1,\ldots,x_r)$ can be written as the summation over all the subsets of $\{1,\dots,r\}$
with $(r-s)$ elements:

\beq\label{eq:11}\ber
\tX^{s}_{n,r}(x_1,\ldots,x_r)=\suml_{\{{i_1},\ldots,{i_{r-s}}\}\sset\{1,\ldots,r\}}x_{i_1}\cdots x_{i_{r-s}}
\tX_{n-r+s,r-s}(x_{i_1},\ldots,x_{i_{r-s}})\\ \ \\
\hspace{2.6cm}=\suml_{\{{i_1},\ldots,{i_{r+1-s}}\}\sset\{1,\ldots,r\}}\frac{1}{s}
\suml_{j\in\{{i_1},\ldots,{i_{r+1-s}}\}}
\frac{x_{i_1}\cdots x_{i_{r+1-s}}}{x_j}\cdot
\tX_{n-r+s,r-s}(x_{i_1},\ldots,\widehat{x_j},\ldots,x_{i_{r+1-s}}).
\eer\eeq
Here in the second line the notation $\widehat{x_j}$ means that the variable $x_j$ is omitted. Note that
in the second line the
summation is as the summation in $\tX^{s-1}_{n,r}(x_1,\ldots,x_r)$, hence these terms can be combined.
From (\ref{eq:11}) one gets
$$X^{s}_{n,r}(x_1,\ldots,x_r)=\suml_{\{{i_1},\ldots,{i_{r+1-s}}\}\sset\{1,\ldots,r\}}\
\suml_{j\in\{{i_1},\ldots,{i_{r+1-s}}\}}
\frac{x_{i_1}\cdots x_{i_{r+1-s}}}{x_j}\cdot
\frac{\tX_{n-r+s,r-s}(x_{i_1},\ldots,\widehat{x_j},\ldots,x_{i_{r+1-s}})}{s\binom{r}{s}\binom{n-1}{r-s-1}}, $$
and by similar argument  $X^{1}_{n,r+1-s}(x_{i_1},\ldots,x_{i_{r+1-s}})$ equals
$$\suml_{\{{i_1},\ldots,{i_{r+1-s}}\}\sset\{{i_1},\ldots,{i_{r+1-s}}\}}\
\suml_{j\in\{{i_1},\ldots,{i_{r+1-s}}\}}
\frac{x_{i_1}\cdots x_{i_{r+1-s}}}{x_j}\cdot
\frac{\tX_{n-r+s,r-s}(x_{i_1},\ldots,\widehat{x_j},\ldots,x_{i_{r+1-s}})}{\binom{r+1-s}{1}\binom{n-1}{r-s-1}}. $$
These two identities  and $s\binom{r}{s}\binom{n-1}{r-s-1}=
\binom{r+1-s}{1}\binom{n-1}{r-s-1}\binom{r}{s-1}$
provide
$$X^{s}_{n,r}(x_1,\ldots,x_r)=\suml_{\{{i_1},\ldots,{i_{r+1-s}}\}\sset\{1,\ldots,r\}}\
\frac{1}{\binom{r}{s-1}}\cdot X^{1}_{n,r+1-s}(x_{i_1},\ldots,x_{i_{r+1-s}}).$$
Using the first line of (\ref{eq:11}), by similar comparison we get
$$X^{s-1}_{n,r}(x_1,\ldots,x_r)=\suml_{\{{i_1},\ldots,{i_{r+1-s}}\}\sset\{1,\ldots,r\}}\
\frac{1}{\binom{r}{s-1}}\cdot X^{0}_{n,r+1-s}(x_{i_1},\ldots,x_{i_{r+1-s}}).$$
Therefore,
$$
X^{s}_{n,r}-X^{s-1}_{n,r}
=\suml_{\{{i_1},\dots,{i_{r+1-s}}\}\sset\{1,\dots,r\}}\frac{1}{\bin{r}{s-1}}
\Big(X^{1}_{n,r+1-s}(x_{i_1},\dots,x_{i_{r+1-s}})-X^{0}_{n,r+1-s}(x_{i_1},\dots,x_{i_{r+1-s}})\Big).
$$
Thus, if the inequality $X^{1}_{n',r'}\ge X^{0}_{n',r'}$ is satisfied for any $n'\le n$ and $r'\le r$,
then $X^{s}_{n',r'}\ge X^{s-1}_{n',r'}$ is also satisfied for any $n'\le n$, $r'\le r$ and  $0\le s\le r'-1$.

{\bf Step 2.} We prove $X^{1}_{n,r}\ge X^{0}_{n,r}$, or, equivalently,
$\frac{\tX^{1}_{n,r}}{r(r-1)}\ge\frac{\tX^{0}_{n,r}}{n-r+1}$.

Note that $\tX^{0}_{n,r}=(\prod^r_{i=1} x_i)\tX_{n-r,r}$ and similarly
$\tX^{1}_{n,r}=(\prod^r_{i=1} x_i)\sum^r_{j=1}\frac{\tX_{n-r+1,r-1}(x_1,\dots,\widehat{x_j},\dots,x_r)}{x_j}$.
Both $\tX_{n-r,r}$ and $\tX_{n-r+1,r-1}$ can be decomposed further according to the $s$--types:
$\tX_{n-r,r}=\sum_{s=0}^{r-1}\tX_{n-r,r}^s$ and
$\tX_{n-r+1,r-1}=\sum^{r-2}_{s=0} \tX^s_{n-r+1,r-1}$. We set
$$\tX^{0,s}_{n,r}:=(\prod^r_{i=1} x_i)\tX_{n-r,r}^s \ \ \mbox{and} \ \
\tX^{1,s}_{n,r}:=(\prod^r_{i=1} x_i)\sum^r_{j=1}\frac{\tX_{n-r+1,r-1}^s(x_1,\dots,\widehat{x_j},\dots,x_r)}{x_j}.$$
We claim that for any  $0\le s\le r-2$ one has:
\beq\label{Eq.6}
(r-s-1)\cdot\tX^{1,s}_{n,r}\ge
(r-s-1)(s+1)\cdot\tX^{0,s+1}_{n,r}+(s+1)(s+2)\cdot\tX^{0,s+2}_{n,r}.
\eeq
This follows from the `elementary' inequality $(I_{i_1i_2}^k)\ : \ x_{i_1}^k+x_{i_2}^k\geq x_{i_1}^{k-1}x_{i_2}+
x_{i_1}x_{i_2}^{k-1}$, where $\{i_1,i_2\}\subset \{1,\ldots,r\}$ and $k\geq 2$.
Indeed, for any fixed pair $(i_1,i_2)$ consider all the monomials of type
$M=\prod_{i=1}^rx_i^{m_i}/(x_{i_1}^{m_{i_1}}x_{i_2}^{m_{i_2}})$ with $m_i> 0$ and of degree $n-k$.
Then $(x_{i_1}^k+x_{i_2}^k)\cdot M\in \tX^{1,s}_{n,r}$. Moreover, each monomial $\ux^{\uk}\in \tX^{1,s}_{n,r}$
can be realized in exactly $(r-1-s)$ ways, where $(r-1-s)$ stays for the number of $k_i$'s with $k_i\geq 2$.

Consider next the same monomial $M$ as before. Then
 $\overline{M}:=(x_{i_1}^{k-1}x_{i_2}+x_{i_1}x_{i_2}^{k-1})\cdot M\in \tX^0_{n,r}$.
If $k>2$ then the number of exponents in $\overline{M}$ which equal 1 is $s+1$, hence
$\overline{M}\in  \tX^{0,s+1}_{n,r}$. If $k=2$ then $\overline{M}\in  \tX^{0,s+2}_{n,r}$.

Any monomial $\ux^{\uk}\in \tX^{0,s+1}_{n,r}$
can be realized in exactly $(r-1-s)(s+1)$ ways, which is $|\{i:k_i\geq 2\}|\cdot
|\{i:k_i=1\}|$, the number of possible candidates for the pair $(i_1,i_2)$ for $(I_{i_1i_2}^k)$.
Furthermore, any monomial $\ux^{\uk}\in \tX^{0,s+2}_{n,r}$
can be realized in exactly $(s+1)(s+2)$ ways, the possible ordered pairs of $\{i:k_i=1\}$.

Hence all the possible inequalities $(I_{i_1i_2}^k)$ multiplied by the possible monomials provide exactly
(\ref{Eq.6}). Note that the above monomial counting is compatible with the identity obtained from $|\cK^s_{n,r}|=\binom{r}{s}\binom{n-1}{r-s-1}$ and
$$r\cdot |\cK^s_{n-r+1,r-1}|=(s+1)\cdot |\cK^{s+1}_{n-r,r}|+
\frac{(s+1)(s+2)}{r-1-s}\cdot |\cK^{s+2}_{n-r,r}|.$$
Now, by taking the sum over $s$ in (\ref{Eq.6}), and by regrouping the right hand side we obtain
$$
\tX^{1}_{n,r}\ge(\prod^r_{i=1} x_i)\sum^{r-2}_{s=0}
\Bigg((s+1)\tX^{s+1}_{n-r,r}+\frac{(s+2)(s+1)}{r-s-1}\tX^{s+2}_{n-r,r}\Bigg)=
(\prod^r_{i=1} x_i)\sum^{r-1}_{s=1}\frac{sr}{r-s+1}\tX^{s}_{n-r,r}
$$
Note that in the last sum the term with $s=0$ can also be
included, as its coefficient vanishes. Thus (for some $c>0$):
\beq\label{Eq.Final.Bound}
c(X^{1}_{n,r}-X^{0}_{n,r})=\frac{n-r+1}{r(r-1)}\cdot
\tX^{1}_{n,r}-\tX^{0}_{n,r}\ge(\prod^r_{i=1} x_i)
\sum^{r-1}_{s=0}\Big(\frac{s(n-r+1)}{(r-s+1)(r-1)}-1\Big)\bin{r}{s}\bin{n-r-1}{r-s-1}\cdot
X^s_{n-r,r}. \eeq Next, the right hand side of
(\ref{Eq.Final.Bound}) is non-negative by the following
generalization of the Chebyshev's sum inequality (which basically
is the summation
$\sum_{s,t}(\alpha_s\alpha_t\beta_s-\alpha_s\alpha_t\beta_t)(x_s-x_t)\geq
0$), see \cite[p. 43]{Hardy-Littlewood-Pólya},
\beq\label{eq:C<gen}
\big(\sum_s \alpha_s\beta_s\big)\big(\sum_s \alpha_s x_s\big)\leq \big( \sum_s\alpha_s\big)
\big(\sum _s\alpha_s\beta_s x_s\big)
\eeq
whenever $x_s$ and $\beta_s$ are both decreasing (or both increasing) sequences and $\alpha_s>0$.
In the present case take $\alpha_s:=\bin{r}{s}\bin{n-r-1}{r-s-1}$, $\beta_s:=\frac{s(n-r+1)}{(r-s+1)(r-1)}$
and $x_s:=X^s_{n-r,r}$.  Clearly $\beta_s$ is increasing, $x_s$ is increasing by induction, and
$\sum_s\alpha_s=\sum_s\beta_s\alpha_s$ by a computation based on
$\sum_s\binom{p}{s}\binom{q}{m-s}=\binom{p+q}{m}$.
Hence, via (\ref{Eq.Final.Bound}), $X^{1}_{n,r}\geq X^{0}_{n,r}$.

\vspace{2mm}

\noindent
{\bf (b)} \  As in part (a), we first reduce the general statement to
$Y^{0}_{n,r,\ell}> Y^{1}_{n,r,\ell}$, and then we prove this particular case too. As in the previous case, set
$\tY^s_{n,r,\l}:=\sum_{\uk\in\cK^s_{n,r}}y_{\uk,\ell}
=\bin{r}{s}\bin{n-1}{r-s-1}Y^s_{n,r,\l}$. We wish to prove:
\beq\label{eq:300}
\frac{\tY^s_{n,r,\l}}{(s+1)(n-r+s+1)}>\frac{\tY^{s+1}_{n,r,\l}}{(r-s)(r-s-1)}.
\eeq

{\bf Step 1.}  Use the decomposition $\cK^s_{n,r}=\coprod \cK^0_{n,r-s}$, the disjoint union
of $\bin{r}{s}$ copies,
to get:
\beq\label{eq:200}
\tY^s_{n,r,\l}=\suml_{\sum k_i=n\atop |\{i:k_i=0\}|=s}\prod^r_{i=1}\frac{1}{(k_i+\l)!}=
\bin{r}{s}\frac{1}{(\l!)^s}\tY^0_{n,r-s,\l}.
\eeq
Similarly, $\tY^{s+1}_{n,r,\l}=\bin{r}{s+1}\frac{1}{(\l!)^{s+1}}\tY^0_{n,r-s-1,\l}=
\bin{r}{s}\frac{1}{(\l!)^{s}}\frac{1}{s+1}\tY^1_{n,r-s,\l}$.
 Here we used $\tY^1_{n,r-s,\l}=\frac{r-s}{\l!}\tY^0_{n,r-s-1,\l}$, cf. (\ref{eq:200}).

Therefore, the inequality (\ref{eq:300}) for any $s$ is equivalent to (\ref{eq:300}) for $s=0$.

\vspace{2mm}

{\bf Step 2.} Here we  prove $\frac{\tY^0_{n,r,\l}}{n-r+1}>\frac{\tY^1_{n,r,\l}}{r(r-1)}$.
We run induction on $n$: we
assume  the stated inequalities, indexed by $(n,r,\ell)$, is true for any $(n',r,\ell')$ with $n'<n$
(but $\ell'$ can be larger than $\ell$). In fact, we use $(n-r,r,\ell+1)\Rightarrow (n,r,\ell)$.

We consider exactly the same combinatorial set--decomposition as in Step 2 of (a), the only difference is that
we replace the inequality $(I_{i_1i_2}^k)$, written for $k\geq 2$,  by $$(I_{i_1i_2,\ell}^k):\ \ \
\frac{2}{(k+\l)!\l!}<\frac{2}{(k-1+\l)!(\l+1)!}.$$
After a computation, we obtain the analogue of (\ref{Eq.Final.Bound}) (for the same positive constant $c$), namely
\beq\label{Eq.Final.Bound.b}
c(Y^{1}_{n,r,\ell}-Y^{0}_{n,r,\ell})=\frac{n-r+1}{r(r-1)}\cdot \tY^{1}_{n,r,\ell}-\tY^{0}_{n,r,\ell}
<\sum^{r-1}_{s=0}\Big(\frac{s(n-r+1)}{(r-s+1)(r-1)}-1\Big)\bin{r}{s}\bin{n-r-1}{r-s-1}\cdot Y^s_{n-r,r,\ell+1}.
\eeq
The right hand side of (\ref{Eq.Final.Bound.b}) is non-positive by (\ref{eq:C<gen}) with reversed inequality,
valid for $\alpha_s>0$, $\beta_s$ and $x_s$ oppositely ordered. Indeed,
 $\alpha_s$ and $\beta_s$ are the same as before with $\beta_s$ is increasing, while $x_s:=Y^s_{n-r,r,\ell+1}$ is a decreasing by induction.
Hence, via (\ref{Eq.Final.Bound.b}), $Y^{1}_{n,r,\ell}< Y^{0}_{n,r,\ell}$.

\vspace{2mm}

{\bf (c)} \ The proof is double induction
over $r$ and $n$. We assume that for any fixed $r$ and $n$ the inequality $(I_{n,r',\ell})$ is true
for any $n$ and $\ell$ and $r'<r$, and  $(I_{n',r,\ell})$ is true for
any  $n'<n$ and any $\ell$. We wish to prove $(I_{n,r,\ell})$.

First we write the left hand side of the inequality as a sum
$$
\sum_{\uk\in\cK_{n,r}} \uy_{\uk,\ell}\cdot \ux^{\uk}\,=\,
\sum_{s=0}^{r-1}\ \sum_{\uk\in\cK_{n,r}^s} \uy_{\uk,\ell}\cdot \ux^{\uk}\ .
$$
Note that corresponding to $s=0$, after we factor out
$x_1\cdots x_r$, the sum over $\cK_{n,r}^0$  can be identified with the left hand side of the inequality
$(I_{n-r,r,\ell+1})$ (multiplied by $x_1\cdots x_r$).
Hence, by the inductive assumption,
$$  \sum_{\uk\in\cK_{n,r}^0} \uy_{\uk,\ell}\cdot \ux^{\uk}\leq Y^0_{n,r,\ell}\cdot
\sum_{\uk\in\cK_{n,r}^0}  \ux^{\uk}.
$$
For $s=1$, the sum over $\cK_{n,r}^1$ is a sum of $r$ sums corresponding to the `missing'
coordinate $x_i$, and each of them can be identified (after factorization of a monomial) with the
inequality $(I_{n-(r-1),r-1,\ell+1})$. For an arbitrary $s\leq r-2$ one can apply in the similar way the
inequality $(I_{n-(r-s),r-s,\ell+1})$. In the case of $s=r-1$ all coefficients
$\uy_{\uk,\ell}$ equal $[(l')^{r-1}(n+l)!]^{-1}$. Therefore, by induction, we get
$$
\sum_{s=0}^{r-1}\,\sum_{\uk\in\cK_{n,r}^s} \uy_{\uk,\ell}\cdot \ux^{\uk}\,\leq\,
\sum_{s=0}^{r-1}|\cK_{n,r}^s|\cdot Y^s_{n,r,\ell}\cdot X^s_{n,r}.
$$
But, using parts (a) and (b), by Chebyshev's sum inequality ((\ref{eq:C<gen}) with $\alpha_s=1$):
$$
 \sum_{s=0}^{r-1}|\cK_{n,r}^s|\cdot Y^s_{n,r,\ell}\cdot X^s_{n,r}\leq
|\cK_{n,r}|\cdot Y_{n,r,\ell}\cdot  X_{n,r},
$$
whose right hand side  is the left hand side of Combinatorial Inequality. This ends the proof of (c).
\epr

\vspace{2mm}

The above discussion and the statement of Theorem \ref{the:2}(b) imply the inequality
(\ref{eq:C<}) from the introduction as well. This is a proof of  (\ref{eq:C<})
in the spirit of the Combinatorial Inequality (based on Chebyshev's type inequalities), for a different proof see
\cite{Kerner-Nemethi.a}.

\bcor\label{ss:ineq} $C_{n,r}>C_{n,r+1}$ for any $n\geq 1$ and $r\geq 1$.
\ecor
\bpr By (\ref{eq:1}) the inequality $C_{n,r}>C_{n,r+1}$ is equivalent to
$Y_{n,r}<Y_{n,r+1}$. We drop the index $\ell=1$ from the notations
(hence, we write e.g.  $Y_{n,r}:=Y_{n,r,1}$), and we set
$Y^{\geq 1}_{n,r}:=\cup_{s\geq 1}Y^s_{n,r}$, and similarly  $\tY^{\geq 1}_{n,r}$ and
$\cK^{\geq 1}_{n,r}$. By \ref{the:2}, part (b), the mean $Y^0_{n,r+1}$ is the largest among
$\{Y^s_{n,r+1}\}_s$, hence $Y_{n,r+1}>Y^{\geq 1}_{n,r}$. Hence, we need to prove $Y^{\geq 1}_{n,r}\geq
Y_{n,r}$. They can be decomposed as sums over the same index set. Indeed, by similar arguments as in the previous proof
part (b)
$$\tY^{\geq 1}_{n,r+1}=\sum _{s=0}^{r-1}\binom{r+1}{s+1}\cdot \tY^0_{n,r-s}=
\sum _{s=0}^{r-1}\frac{\binom{r+1}{s+1}}{\binom{r}{s}}\cdot\tY^s_{n,r}=
\sum _{s=0}^{r-1}\frac{\binom{r+1}{s+1}}{\binom{r}{s}}\cdot\binom{r}{s}\binom{n-1}{r-s-1}\cdot Y^s_{n,r}.
$$
This can be rewritten as
\beq \label{eq:c<1}
Y^{\geq 1}_{n,r+1}=\frac{1}{|\cK^{\geq 1}_{n,r+1}|}\cdot
\sum _{s=0}^{r-1}\frac{r+1}{s+1}\binom{r}{s}\binom{n-1}{r-s-1}\cdot Y^s_{n,r}.
\eeq
Similarly,
\beq \label{eq:c<2}
Y_{n,r}=\frac{1}{|\cK_{n,r}|}\cdot
\sum _{s=0}^{r-1}\binom{r}{s}\binom{n-1}{r-s-1}\cdot Y^s_{n,r}.
\eeq
Hence, via (\ref{eq:c<1}) and (\ref{eq:c<2}),
the inequality $Y_{n,r}\leq Y^{\geq 1}_{n,r}$ follows from (\ref{eq:C<gen}) applied for  $\alpha_s=\binom{r}{s}\binom{n-1}{r-s-1}$, $\beta_s=\frac{r+1}{s+1}$ and $x_s=Y^s_{n,r}$.
Indeed,   $\alpha_s>0$ while   $x_s$ and $\beta_s$ are both decreasing sequences.
For $x_s$ this  follows from Theorem \ref{the:2}(b).
\epr

\section{The homogeneous case for $n=2$}\label{s:4}

Assume that $(X,0)$ is a 2--dimensional ICIS with  $p_i\geq 2$ for all $i$.
Then $C_{2,r}=4\frac{r+1}{r+\frac{1}{3}}$.  Set
\[P:=\prod_ip_i,\ \text{ the multiplicity of }(X,0).\]
 Then using \ref{Thm.Milnor.Genus.For.ICIS} one obtains
$$\frac{p_g}{P}= \sum_i\frac{(p_i-1)(p_i-2)}{6}+\sum_{i<j}\frac{(p_i-1)(p_j-1)}{4}$$
 and $$\frac{\mu+1-P}{P}=\sum_i\big((1-p_i)+(1-p_i)^2\big)+\sum_{i<j}(1-p_i)(1-p_j).$$
By a computation
$\mu+P\cdot E +1= C_{2,r}\cdot p_g$, where $E:=
\frac{r-1}{3r+1}\sum^r_{i=1}(p_i-1)-\sum_{i<j}\frac{(p_i-p_j)^2}{3r+1}-1$.

If $r=1$ then $E=-1$, but for $r\geq 2$ and for some choices of $p_i$'s
(e.g. whenever they are all equal) $E$ might be positive, providing
$C_{2,r}\cdot p_g\geq \mu+1$. We collect here the precise statements:

\bthe (a) \ If $r=1$ then $6p_g = \mu+1-P$.

(b) If $r\geq 2$ then  $C_{2,r}\cdot  p_g=4\cdot \frac{r+1}{r+1/3}\cdot p_g$.
In general the bound  $C_{2,r}\cdot  p_g\leq \mu+1$ does \underline{not} hold, although asymptotically
$\frac{\mu}{p_g}$ tends to $\frac{4(r+1)}{r+1/3}$.

(c) For any $r$ the inequality $4p_g\leq \mu+1-P$ is valid, and if $p_i=2$ is allowed then
4 is the sharpest bound whenever $r>1$.

(d) If $p_i\geq d+1$ for all $i$ then
$$4\cdot\frac{d(r-1)+2(d-1)} {d (r-1)+\frac{4}{3}(d-1)}\cdot p_g\leq \mu+1-P.$$
\ethe
\bpr Using the above explicit formulae, all the statements are elementary.
Let us give some hint for (d). Note that
$(4+c)p_g\leq \mu+1-P$  reads as
\begin{equation}\label{n8}
(4+c)\sum_i\frac{(p_i-1)(p_i-2)}{6}+c\sum_{i<j}\frac{(p_i-1)(p_j-1)}{4}\leq
\sum_i(1-p_i)+\sum_{i}(1-p_i)^2.\end{equation}
Using
$\sum_{i=1}^ra_i^2\geq \frac{2}{r-1}\sum_{i<j}a_ia_j$, the inequality
(\ref{n8}) follows from
$$(4+c)\sum_i\frac{(p_i-1)(p_i-2)}{6}+c\cdot \frac{r-1}{8}\sum_i(p_i-1)^2\leq
\sum_i(1-p_i)+\sum_{i}(1-p_i)^2.$$
This is a sum over $i$ of elementary quadratic inequalities
whose discussion is left to the reader.
\epr


\begin{thebibliography}{99}
\bibitem[Abramowitz-Stegun]{Abramowitz-Stegun}
M. Abramowitz, I.A. Stegun, {\em Handbook of mathematical functions with formulas, graphs, and mathematical tables.}
 National Bureau of Standards Applied Mathematics Series, 55, 1964.


\bibitem[AGLV-book]{AGLV} V.I. Arnol'd, V.V. Goryunov, O.V. Lyashko, V.A. Vasil'ev, {\em
Singularity theory.I.} Reprint of the original English edition
from the series Encyclopaedia of Mathematical Sciences [Dynamical
systems. VI, Encyclopaedia Math. Sci., 6, Springer, Berlin, 1993].
Springer-Verlag, Berlin, 1998.


\bibitem[Ashikaga1992]{Ashikaga92}T. Ashikaga, {\em Normal two-dimensional hypersurface triple points
and the Horikawa type resolution.} Tohoku Math. J. (2) 44 (1992), no. 2, 177--200.

\bibitem[Buchweitz-Greuel1980]{Buchweitz-Greuel-1980}R.-O. Buchweitz, G.-M. Greuel, {\em
The Milnor number and deformations of complex curve singularities.} Invent. Math. 58 (1980), no. 3, 241--281.


\bibitem[Durfee1978]{Durfee1978} A.H. Durfee, {\it The signature of smoothings of complex surface singularities.}
Math. Ann. 232 (1978), no. 1, 85--98.


\bibitem[Greuel1975]{Greuel1975}G.M. Greuel, {\em Der Gauss-Manin-Zusammenhang isolierter Singularit\"{a}ten von
vollst\"{a}ndigen Durchschnitten.} Math. Ann. 214 (1975), 235-266.

\bibitem[Greuel-Hamm1978]{Greuel-Hamm1978}G.M. Greuel, H.A. Hamm, {\em Invarianten quasihomogener vollst\"{a}ndiger
Durchschnitte.}, Invent. Math. 49 (1978), no. 1, 67-86.


\bibitem[Hamm1986]{Hamm1986}H.A. Hamm, {\em Invariants of weighted homogeneous singularities.} Journ\'{e}es Complexes 85
(Nancy, 1985), 613, Inst. \'{E}lie Cartan, 10, Univ. Nancy, Nancy, 1986.

\bibitem[Hamm2011]{Hamm2011}H.A. Hamm, {\em Differential forms and Hodge numbers for toric complete intersections},
 arXiv:1106.1826.

\bibitem[Hardy-Littlewood-P\'{o}ólya]{Hardy-Littlewood-Pólya}G.H.Hardy, J.E.Littlewood, G.P\'{o}lya,
{\em Inequalities}. Reprint of the 1952 edition. Cambridge Mathematical Library.
Cambridge University Press, Cambridge, 1988.


\bibitem[Jordan1965]{Jordan1965}Ch.Jordan, {\em Calculus of finite differences}. Third Edition.
Introduction by Harry C. Carver Chelsea Publishing Co., New York 1965 xxi+655 pp

\bibitem[Kerner-N\'emethi2009]{Kerner-Nemethi2009} D. Kerner and A. N\'emethi,
{\em The Milnor fibre signature is not semi-continuous},
Proc. of the Conference in Honor of the 60th Birthday of A. Libgober,
Topology of Algebraic Varieties, Jaca (Spain),  June 2009;
 Contemporary Math. 538 (2011), 369--376.

\bibitem[Kerner-N\'emethi.a]{Kerner-Nemethi.a} D. Kerner and A. N\'emethi,
{\em A counterexample to Durfee conjecture}, Comptes Rendus Math\'{e}ématiques
 de l'Acad\'{e}émie des Sciences, vol.34 (2012), no.2.  arXiv:1109.4869



\bibitem[Kerner-N\'emethi.c]{Kerner-Nemethi.c} D. Kerner, A. N\'emethi,
{\em On Milnor number and singularity genus}, in preparation.

\bibitem[Khovanskii1978]{Khovanskii1978} A.G. Khovanskii,
{\em Newton polyhedra, and the genus of complete intersections.}
 (Russian) Funktsional. Anal. i Prilozhen. 12 (1978), no. 1, 51-61.

\bibitem[Laufer1977]{Laufer1977} H.B. Laufer, {\em  On $\mu$ for surface singularities,}
 Proceedings of Symposia in Pure Math.  30, 45-49,  1977.

\bibitem[Looijenga-book]{Looijenga-book} E. Looijenga, {\em Isolated Singular Points on Complete Intersections.}
 London Math. Soc. LNS 77, CUP, 1984.

\bibitem[Looijenga1986]{Looijenga1986} E. Looijenga, {\em Riemann-Roch and smoothings of singularities.}
 Topology 25 (1986), no. 3, 293--302.


\bibitem[Morales1985]{Morales1985}M. Morales, {\em Fonctions de Hilbert, genre g\'{e}om\'{e}trique
d'une singularit\'{e} quasi homog\`{e}ne Cohen-Macaulay.} C. R. Acad. Sci. Paris S\'{e}r. I Math. 301 (1985),
no. 14, 699--702.


\bibitem[Melle-Hern\'{a}ndez2000]{Melle00}A. Melle-Hern\'{a}ndez, {\it Milnor numbers for surface singularities.}
  Israel J. Math.  115  (2000), 29--50.

\bibitem[Milnor-book]{Milnor-book} J. Milnor, {\em Singular points of complex hypersurfaces,}
 Annals of Math. Studies  61, Princeton University Press 1968.

\bibitem[N\'emethi98]{Nemethi98} A. N\'emethi, {\em Dedekind sums and the signature of $f(x,y)+z^N$}.
Selecta Math. (N.S.) 4 (1998), no. 2, 361--376.

\bibitem[N\'emethi99]{Nemethi99} A. N\'emethi,    {\em  Dedekind sums and the signature of $f(x,y)+z^N$,II.}
Selecta Math. (N.S.)  5 (1999), 161--179.


\bibitem[Saito1981]{Saito1981}M. Saito, {\em On the exponents and the geometric genus of an isolated hypersurface
singularity}. Singularities, Part 2 (Arcata, Calif., 1981), 465--472,
Proc. Sympos. Pure Math., 40, Amer. Math. Soc., Providence, RI, 1983.

\bibitem[Seade-book]{Seade}J. Seade, {\it On the Topology of Isolated Singularities in Analytic Spaces.} Progress in
Mathematics 241, Birkh\"auser 2006.


\bibitem[Tomari1993]{Tomari93} M. Tomari, {\em The inequality $8p_g<\mu$ for hypersurface two-dimensional
isolated double points.} Math. Nachr. 164 (1993), 37--48.


\bibitem[Wahl1981]{Wahl 1981}J. Wahl, {\em Smoothings of normal surface singularities}, Topology 20 (1981), 219--246.

\bibitem[Xu-Yau1993]{Xu-Yau1993}Y.-J. Xu, S.S.-T. Yau, {\em Durfee conjecture and coordinate free characterization of
homogeneous singularities.} J. Differential Geom. 37 (1993), no. 2, 375--396.

\bibitem[Yau-Zhang2006]{Yau-Zhang2006}St.-T.Yau, L.Zhang, {\em An upper estimate of integral points in real simplices with
an application to singularity theory.} Math. Res. Lett. 13 (2006), no. 5--6, 911--921.
\end{thebibliography}
\end{document}